\documentclass[12pt]{article}
\usepackage{amsmath}
\usepackage{amssymb}
\usepackage{geometry}
\usepackage{amscd}
\usepackage{xypic}
\usepackage{mathrsfs}
\usepackage{amsthm}
\usepackage{color}
\usepackage{fancyhdr}
\usepackage{graphicx}
\usepackage{pgf,tikz}
\usetikzlibrary{arrows}
\usetikzlibrary{calc,intersections,through,backgrounds}
\usepackage{array}
\usepackage{float}
\usepackage{hyperref}

\newtheorem{theorem}{Theorem}[section]
\newtheorem{definition}{Definition}[section]
\newtheorem{proposition}{Proposition}[section]

\theoremstyle{definition}

\newtheorem{remark}{Remark}[section]

\begin{document}

\title{Invariants, Bitangents and Matrix Representations of Plane Quartics with 3-Cyclic Automorphisms 
}

\author{Dun Liang      
}

\maketitle

\begin{abstract} In this work we compute the Dixmier invariants and bitangents of the plane quartics with 3,6 or 9-cyclic automorphisms, we find that a quartic curve with 6-cyclic automorphism will have 3 horizontal bitangents which form an asysgetic triple. We also discuss the linear matrix representation problem of such curves, and find a degree 6 equation  of 1 variable which solves the symbolic solution of the linear matrix representation problem for the curve with 6-cyclic automorphism. 
\end{abstract}
\section{Introduction}

The study of the geometry of plane quartics is one of the most beautiful achievements in classical algebraic geometry. Back to the late $19^{\rm th}$  and early $20^{\rm th}$ century, there were many studies on the existence and configurations of the 28 bitangents of a plane quartic such as \cite{Jacobi}, \cite{Hesse1}, \cite{Hesse2}, \cite{Steiner}, and so on.
For the invariants of plane quartics, Shioda computed the invariants ring in \cite{Shioda}.
 However, the algebraic invariants of plane quartics were found much later by 
\cite{Dixmier}, \cite{Ohno} and \cite{Elsenhans}.
 In this work we compute the invariants and bitangents of plane quartics with 3-cyclic automorphism, and discuss the linear matrix representation problem (see \cite{Helton},\cite{Vinnikov1},\cite{Vinnikov2}) of such curves. The classification of automorphism was given by
\cite{Kantor} and \cite{Wiman}. There are many places one can see the full list, for example, in Section 6.5 of \cite{CAG}. 

Explicitly, we consider the curves
\begin{align}
&C_3=C_3(r,s): &y^3=x(x-1)(x-r)(x-s) \\
&C_6=C_6(r): &y^3=x(x-1)(x-r)(x-1+r) \\
&C_9:&y^3=x(x^3-1)
\end{align}
with automorphism group ${\mathbb Z}/3$, ${\mathbb Z}/6$ and ${\mathbb Z}/9$ respectively.

We compute the invariants of these curves. The curves $C_6$ and $C_9$ are special cases of $C_3$. Thus we also compute the cutting equations of the invariants of $C_6$ and $C_9$ as special cases of $C_3$. In modern point of view, a smooth plane quartic is the canonical model of a smooth projective non-hyperelliptic curve of genus $3$. Let ${\cal M}_3$ be the moduli space of projective curves of genus $3$, and let ${\cal M}_3^{\rm non}$ be the non-hyperelliptic locus of ${\cal M}_3$. Then the Dixmier invariants $I_3,I_6,I_9,I_{12},I_{15}, I_{18}, I_{27}$, as functions of the coefficients of a given ternary quartic like an analog of the $j$-invariant of a given cubic curve, also could be regarded as the coordinates of ${\cal M}_3^{\rm non}$. Let $G$ be a finite group. If we write $X^G$ as the subvariety of ${\cal M}_3^{\rm non}$ parametrizing curves with automorphism group containing $X^G$, then we have $X^{{\mathbb Z}/9}\subset X^{{\mathbb Z}/3}\subset {\cal M}_3^{\rm non}$ and $X^{{\mathbb Z}/6}\subset X^{{\mathbb Z}/3}\subset {\cal M}_3^{\rm non}$. In this point of view, we are trying to find the ``defining equation" of $X^{{\mathbb Z}/3}$, $X^{{\mathbb Z}/6}$ and $X^{{\mathbb Z}/9}$ in ${\cal M}_3^{\rm non}$. 

The explicit formulae of the Dixmier invariants are listed in Section \ref{Dix}.
We use {\tt Maxima} to compute the Dixmier invariants. We will give the {\tt wxMaxima} notebook in the appendix.
Summarizing Section \ref{C369Dix}, we have the following
\begin{theorem}
The curve $C_3$ satisfies that 
$$I_3=I_6=I_{12}=I_{15}=0$$
and $I_9,I_{18}$ are algebraically independent. For the curve $C_6$, the invariants $I_9$ and $I_{18}$ satisfies a degree 8 affine equation. Furthermore, the curve $C_9$ is the curve that all Dixmier invariants vanish.
\end{theorem}

The algebraic conditions between the invariants of $C_6$ are computed by {\tt Macaulay2} \cite{Macaulay2}.

We use the idea in \cite{Sturmfels} to compute the bitangents of a plane quartic. This program is also realized by {\tt Macaulay2}. We summarize Section \ref{c369bit} as the following theorem.
\begin{theorem}
The curve $C_9$ has all 28 explicit symbolic solutions for the bitangents. The curve $C_6$ has 3 horizontal explicit bitangents which form a triple of asysgetic sets. 
\end{theorem}

The definition of asysgetic sets come from the theory of theta characteristics. For details one can see \cite{Mumford}, our definition in Section \ref{bitan} is a geometric description as in \cite{Sturmfels}.

For the linear matrix determinant representation problem of such curves, we use the idea in \cite{Sturmfels2}. The problem asks wether the equation of a plane curve $C$ could be written of the form
$$\det(xA+yB+zC)$$ for some symmetric matrices $A,B,C$ of constants. Our result is Theorem \ref{5.1} as follows:
\begin{theorem}
The matrix representation problem for $C_6$ has a symbolic solution as a solution of a degree 6 polynomial of 1 variable whose coefficients are defined over $\overline{K(r,s)}$.
\end{theorem}

\section{Automorphisms of Plane Quartics}

We consider the algebraic varieties over the algebraic closure $K=\overline{\mathbb Q}$ of the rational field $\mathbb Q$ in the complex number $\mathbb C$ since we are interested in the geometric properties of such varieties.
 However, some of the algorithms we use later in this work will be realized over $\mathbb Q$ only. 
 In this section, let $K=\overline{\mathbb Q}$.

Let $C$ be a smooth projective curve over $K$. If the genus $g(C)$ of $C$ is 3, and $C$ is non-hyperelliptic, then the canonical model of $C$ is a plane quartic and is isomorphic to $C$. 
Let $x,y,z$ be the coordinates of the projective plane ${\mathbb P}^2$, if we want to emphasize the coordinates, we also write ${\mathbb P}^2$ as ${\mathbb P}^2_{(x,y,z)}$.
 Let $k[x,y,z]_d$ be the homogeneous degree $d$-part of the polynomial ring $K[x,y,z]$.
  Thus $k[x,y,z]_d\simeq {\rm Sym}^d((K^\vee)^3)$, the 3rd symmetric product of $K^\vee = {\rm Hom}_K(K,K)\simeq K$.
   We write ${\cal P}_n^d:={\rm Sym}^d((K^\vee)^n)$.
    Thus, let $F_C=F_c(x,y,z)$ be the equation of $C$, we say both $F_C\in {\cal P}_3^4$ and $F_C\in K[x,y,z]_4$.
    
    An element $F\in K[x,y,z]_4$ should be written as
    $$F(x,y,z)=\sum_{i+j+k=4}a_{ijk}x^iy^jz^k.$$
    
    Let $C$, $D$ be two smooth non-hyperelliptic genus $g$ curves over $K$.
    The canonical models $\kappa_C$, $\kappa_D$ of $C$ and $D$ are closed subvarieties of degree $2g-2$ in ${\mathbb P}^{g-1}$. 
    Since $C$ and $D$ are non-hyperelliptic, we have $C\simeq \kappa_C$ and $D\simeq \kappa_D$. 
    The theory of algebraic curves says that $C$ and $D$ are isomorphic as algebraic varieties if and only if $\kappa_C$ could be transformed to $\kappa_D$ by a non-degenerated projective linear transformation on the coordinates of ${\mathbb P}^{g-1}$.
    In particular, an automorphism of a non-hyperelliptic curve $C$ is a projective automorphism on the canonical model $\kappa_C$ of $C$.
    
    In this work we consider non-hyperelliptic genus 3 curves with cyclic automorphism groups ${\mathbb Z}/3$, ${\mathbb Z}/6$ and ${\mathbb Z}/9$.
    
    The genus 3 non-hyperelliptic curves with ${\mathbb Z}/3$-automorphisms is a 2-dimensional family
$$C_3=C_3(r,s):\quad y^3=x(x-1)(x-r)(x-s).$$
This is a family of smooth quartics written on the affine chart $\{z=1\}$ of the projective plane ${\mathbb P}^2_{(x,y,z)}$ with $K$-parameters $r$ and $s$. 

Also we have the 1-dimensional family
$$C_6=C_6(r): \quad y^3=x(x-1)(x-r)(x-1+r)$$ of curves with automorphism group ${\mathbb Z}/6$ and the curve
$$C_9:\quad y^3=x(x^3-1)$$
whose automorphsm group is ${\mathbb Z}/9$

In the following sections we will compute the invariants and bitangents of $C_3$, $C_6$ and $C_9$.

\section{Dixmier Invariants of $C_3$, $C_6$ and $C_9$}
\subsection{Dixmier Invariants of Plane Quartics}\label{Dix}

Our notations follows from \cite{Kohel}. First we introduce some notations. In general, let $f\in K[x_1,\ldots, x_n]$ be a polynomial, we use $D_f$ to denote the differential operator determined by $f$. Explicitly, let 
\begin{equation}\label{genf} f=f(x_1,\ldots ,x_n)= \sum_{(i_1,\ldots , i_n)\in {\mathbb Z}^n_+} a_{i_1,\ldots ,i_n}x_1^{i_1}\cdots x_n^{i_n}\end{equation}
where $a_{i_1,\ldots ,i_n}\in K$ be the coefficient of the monomial $x_1^{i_1}\cdots x_n^{i_n}$ for $(i_1,\ldots , i_n)\in {\mathbb Z}^n_+$ and (\ref{genf}) be a finite sum. For the rest of this paper, we will not emphasize that the powers $i_1,\ldots , i_n$ are non-negative integers again. 

The map $D_f$ means
$$\begin{array}{cccc}
D_f : & K[x_1,\ldots, x_n] & \longrightarrow & K[x_1,\ldots, x_n] \\ \\ 
& g(x_1,\ldots , x_n) & \longmapsto & \displaystyle{\sum_{(i_1,\ldots , i_n)\in {\mathbb Z}^n_+} a_{i_1,\ldots ,i_n}\frac{\partial^{i_1+\cdots + i_n}}{\partial x_1^{i_1}\cdots{\partial x_n^{i_n}}} g(x_1,\ldots , x_n)}.
\end{array}
$$
If we use $D(f,g)$ to denote $D_f(g)$, $\forall f,g\in K[x_1,\ldots ,x_n]$, then the map 
$$D:K[x_1,\ldots ,x_n] \times K[x_1,\ldots ,x_n] \longrightarrow K[x_1,\ldots ,x_n]$$ has some obvious properties as the following:

\begin{itemize}
\item $D$ is bilinear.
\item Let ${\rm deg}(f)$ be the degree of $f$ for all $f\in K[x_1,\ldots ,x_n]$.
 Let $f,g\in K[x_1,\ldots ,x_n]$. If ${\rm deg}(f)> {\rm deg}(g)$, then $D_f(g)=0$. If 
${\rm deg}(f)> {\rm deg}(g)$, then $D_f(g)\leq {\rm deg}(g)-{\rm deg}(f)$.  Let $f=x_1^{i_1}\cdots x_n^{i_n}$ and $g= x_1^{j_1}\cdots x_n^{j_n}$ be two monomials such that ${\rm deg}(f)= {\rm deg}(g)$, then $D_f(g)=i_1!\cdots i_n!\delta_{fg}$ where $\delta_{fg}$ is the Kronecker delta of $f$ and $g$.
\end{itemize}

For any $f\in K[x_1,\ldots ,x_n]$, let $H(f)$ be the half Hessian matrix of $f$. For example, if $f\in K[x,y,z]$, then 
$$H(f)=\begin{pmatrix}
\displaystyle{\frac{\partial^2}{\partial x^2}} & \displaystyle{\frac{\partial^2}{\partial x \partial y}} & \displaystyle{\frac{\partial^2}{\partial x \partial z}} \\ 
 \displaystyle{\frac{\partial^2}{\partial x \partial y}}& \displaystyle{\frac{\partial^2}{\partial y^2}} & \displaystyle{\frac{\partial^2}{\partial y \partial z}} \\
\displaystyle{\frac{\partial^2}{\partial x \partial z}} & \displaystyle{\frac{\partial^2}{\partial y \partial z}} & \displaystyle{\frac{\partial^2}{\partial z^2}}
\end{pmatrix}.$$
Let $H^*(f)$ be the adjoint of $H(f)$.

Another notation is the dot product of two matrices. Let $A=(a_{ij})_{n\times n}$ and $B=(b_{ij})_{n\times n}$ be two $n\times n$ matrices. 
Then the dot product ``$\langle\, ,\rangle$" is defined by 
$$\langle\, A , B \, \rangle := \sum_{1\leq i,j\leq n} a_{ij}b_{ji}.$$

With these notations, we describe the Dixmier invariants of plane quartics.

Let $f,g\in K[x,y,z]_2$ be two quadratic homogeneous polynomials. 
Define
\begin{align*}
&J_{1,1}(f,g)= \langle\, H(f), H(g)\,\rangle , \\
&J_{2,2}(f,g)= \langle\, H^*(f), H^*(g)\,\rangle , \\
&J_{3,0}(f,g) = J_{3,0}(f)= \det (H(f)) ,\\
&J_{0,3}(f,g)=J_{0,3}(g) = \det(H(g)).
\end{align*}

Let $F\in K[x,y]_r$, $G\in K[x,y]_s$ be two homogeneous polynomials of degree $r$ and $s$, respectively. For $k\leq \min\{r,s\}$, define 
\begin{equation}\label{FGK}
(F,G)^k:= \frac{(r-k)!(s-k)!}{r!s!} \left.\left(\frac{\partial^2}{\partial x_1\partial y_2}- \frac{\partial^2}{\partial y_1\partial x_2}\right)^kF(x_1,y_1)G(x_2,y_2)\right|_{(x_i,y_i)=(x,y),i,1,2
}
\end{equation}

Let $P=P(x,y)\in K[x,y]_4$ be a quartic binary form. Let $Q=(P,P)^4$ defined as (\ref{FGK}). Also we let
\begin{equation}\label{SPD}\begin{split}
&\Sigma(P)=\frac{1}{2} (P,P)^4,\quad  \Psi(P)= \frac{1}{6} (P,Q)^4 \\
& \Delta(P) = \Sigma(P)^3 - 27 \Psi(P)^2\end{split}
\end{equation}
Then $\Delta(P)$ is the discriminant of $P$.

Let $u,v$ be two $K$-variables. For quartic $f\in K[x,y,z]_4$, let 
$$g=g(x,y)=f(x,y,-ux-vy).$$
Then $g(x,y)$ is a homogeneous polynomial of degree 4 with respect to the variables $x$ and $y$, and the coefficients of $g$ are expressions of $u$ and $v$. 
Thus we can define $\Sigma(g)$ and $\Psi(g)$ as in (\ref{SPD}). 
Since $\Sigma$ and $\Psi$ are expressions of the coefficients, 
we have $\Sigma(g)$ and $\Psi(g)$ are expressions of $u$ and $v$. 
An explicit computation shows that $\Sigma(g)$ and $\Psi(g)$ are polynomials of degree $2$ and $3$ in the polynomial ring $K[u,v]$ respectively.
Let $\sigma(u,v,w)$ and $\psi(u,v,w)$ be the homogenization of $\Sigma(g)$ for $w$, and $\psi(u,v,w)$ be the homogenization of $\Psi(g)$ for $w$. 
Then $\sigma(u,v,w)\in K[u,v,w]_2$ and $\psi(u,v,w)\in K[u,v,w]_3$.
Finally, we substitute $u=x, v=y, w=z$ into $\sigma(u,v,w)$ and $\psi(u,v,w)$. 
For $f\in K[x,y,z]_4$, we define
\begin{equation}\label{sp}
\begin{split}
\sigma(f)=\sigma=\sigma(x,y,z)\in K[x,y,z]_2 \\
\psi(f)=\psi=\psi(x,y,z)\in K[x,y,z]_3
\end{split}
\end{equation}
\begin{definition} Let $f\in K[x,y,z]_4$, let $\sigma$, $\psi$ defined as in (\ref{sp}). Let $\rho=D_f(\psi)$ and $\tau=D_\rho(f)$.
The Dixmier invariants are defined as 
\begin{equation}
\begin{split}
&I_3=D_\sigma(f), \quad I_9= J_{1,1}(\tau, \rho), \quad I_{15}=J_{3,0}(\tau), \\
& I_6 = D_\psi (H)-8I_3^2, \quad I_{12}=J_{0,3}(\rho), \quad I_{18}= J_{2,2}(\tau, \rho) \\
& I_{27} = \Delta = \sigma^3-27\psi^2
\end{split}
\end{equation}
\end{definition}

\subsection{The Dixmier Invariants of $C_3$, $C_6$ and $C_9$}\label{C369Dix}

We use Maxima to compute the Dixmier invariants of $C_3$, $C_6$ and $C_9$. And we use elimination in Macaulay2 to compute the conditions of the invariants with certain automorphisms.

\begin{proposition} The Dixmier invariants of 
$$C_3(r,s): y^3=x(x-1)(x-r)(x-s)$$
are

$$I_3=I_6=I_{12}=I_{15}=0$$

\begin{flushleft}
$I_9=-\frac{{{r}^{3}}\, {{s}^{5}}}{55296}+\frac{{{r}^{2}}\, {{s}^{5}}}{36864}+\frac{r\, {{s}^{5}}}{36864}-\frac{{{s}^{5}}}{55296}-\frac{77 {{r}^{4}}\, {{s}^{4}}}{331776}+\frac{169 {{r}^{3}}\, {{s}^{4}}}{331776}-\frac{97 {{r}^{2}}\, {{s}^{4}}}{110592}+\frac{169 r\, {{s}^{4}}}{331776}-\frac{77 {{s}^{4}}}{331776}-\frac{{{r}^{5}}\, {{s}^{3}}}{55296}+\frac{169 {{r}^{4}}\, {{s}^{3}}}{331776}+\frac{{{r}^{3}}\, {{s}^{3}}}{13824}+\frac{{{r}^{2}}\, {{s}^{3}}}{13824}+\frac{169 r\, {{s}^{3}}}{331776}-\frac{{{s}^{3}}}{55296}+\frac{{{r}^{5}}\, {{s}^{2}}}{36864}-\frac{97 {{r}^{4}}\, {{s}^{2}}}{110592}+\frac{{{r}^{3}}\, {{s}^{2}}}{13824}-\frac{97 {{r}^{2}}\, {{s}^{2}}}{110592}+\frac{r\, {{s}^{2}}}{36864}+\frac{{{r}^{5}} s}{36864}+\frac{169 {{r}^{4}} s}{331776}+\frac{169 {{r}^{3}} s}{331776}+\frac{{{r}^{2}} s}{36864}-\frac{{{r}^{5}}}{55296}-\frac{77 {{r}^{4}}}{331776}-\frac{{{r}^{3}}}{55296}$

\

$I_{18}=\frac{{{r}^{6}}\, {{s}^{10}}}{402653184}-\frac{{{r}^{5}}\, {{s}^{10}}}{134217728}+\frac{{{r}^{4}}\, {{s}^{10}}}{67108864}-\frac{7 {{r}^{3}}\, {{s}^{10}}}{402653184}+\frac{{{r}^{2}}\, {{s}^{10}}}{67108864}-\frac{r\, {{s}^{10}}}{134217728}+\frac{{{s}^{10}}}{402653184}+\frac{{{r}^{7}}\, {{s}^{9}}}{1358954496}-\frac{163 {{r}^{6}}\, {{s}^{9}}}{10871635968}+\frac{19 {{r}^{5}}\, {{s}^{9}}}{1207959552}-\frac{155 {{r}^{4}}\, {{s}^{9}}}{10871635968}-\frac{155 {{r}^{3}}\, {{s}^{9}}}{10871635968}+\frac{19 {{r}^{2}}\, {{s}^{9}}}{1207959552}-\frac{163 r\, {{s}^{9}}}{10871635968}+\frac{{{s}^{9}}}{1358954496}+\frac{229 {{r}^{8}}\, {{s}^{8}}}{48922361856}-\frac{539 {{r}^{7}}\, {{s}^{8}}}{24461180928}+\frac{2711 {{r}^{6}}\, {{s}^{8}}}{24461180928}-\frac{13241 {{r}^{5}}\, {{s}^{8}}}{97844723712}+\frac{20231 {{r}^{4}}\, {{s}^{8}}}{97844723712}-\frac{13241 {{r}^{3}}\, {{s}^{8}}}{97844723712}+\frac{2711 {{r}^{2}}\, {{s}^{8}}}{24461180928}-\frac{539 r\, {{s}^{8}}}{24461180928}+\frac{229 {{s}^{8}}}{48922361856}+\frac{{{r}^{9}}\, {{s}^{7}}}{1358954496}-\frac{539 {{r}^{8}}\, {{s}^{7}}}{24461180928}+\frac{1913 {{r}^{7}}\, {{s}^{7}}}{48922361856}-\frac{5927 {{r}^{6}}\, {{s}^{7}}}{32614907904}-\frac{1705 {{r}^{5}}\, {{s}^{7}}}{97844723712}-\frac{1705 {{r}^{4}}\, {{s}^{7}}}{97844723712}-\frac{5927 {{r}^{3}}\, {{s}^{7}}}{32614907904}+\frac{1913 {{r}^{2}}\, {{s}^{7}}}{48922361856}-\frac{539 r\, {{s}^{7}}}{24461180928}+\frac{{{s}^{7}}}{1358954496}+\frac{{{r}^{10}}\, {{s}^{6}}}{402653184}-\frac{163 {{r}^{9}}\, {{s}^{6}}}{10871635968}+\frac{2711 {{r}^{8}}\, {{s}^{6}}}{24461180928}-\frac{5927 {{r}^{7}}\, {{s}^{6}}}{32614907904}+\frac{20383 {{r}^{6}}\, {{s}^{6}}}{32614907904}-\frac{35327 {{r}^{5}}\, {{s}^{6}}}{97844723712}+\frac{20383 {{r}^{4}}\, {{s}^{6}}}{32614907904}-\frac{5927 {{r}^{3}}\, {{s}^{6}}}{32614907904}+\frac{2711 {{r}^{2}}\, {{s}^{6}}}{24461180928}-\frac{163 r\, {{s}^{6}}}{10871635968}+\frac{{{s}^{6}}}{402653184}-\frac{{{r}^{10}}\, {{s}^{5}}}{134217728}+\frac{19 {{r}^{9}}\, {{s}^{5}}}{1207959552}-\frac{13241 {{r}^{8}}\, {{s}^{5}}}{97844723712}-\frac{1705 {{r}^{7}}\, {{s}^{5}}}{97844723712}-\frac{35327 {{r}^{6}}\, {{s}^{5}}}{97844723712}-\frac{35327 {{r}^{5}}\, {{s}^{5}}}{97844723712}-\frac{1705 {{r}^{4}}\, {{s}^{5}}}{97844723712}-\frac{13241 {{r}^{3}}\, {{s}^{5}}}{97844723712}+\frac{19 {{r}^{2}}\, {{s}^{5}}}{1207959552}-\frac{r\, {{s}^{5}}}{134217728}+\frac{{{r}^{10}}\, {{s}^{4}}}{67108864}-\frac{155 {{r}^{9}}\, {{s}^{4}}}{10871635968}+\frac{20231 {{r}^{8}}\, {{s}^{4}}}{97844723712}-\frac{1705 {{r}^{7}}\, {{s}^{4}}}{97844723712}+\frac{20383 {{r}^{6}}\, {{s}^{4}}}{32614907904}-\frac{1705 {{r}^{5}}\, {{s}^{4}}}{97844723712}+\frac{20231 {{r}^{4}}\, {{s}^{4}}}{97844723712}-\frac{155 {{r}^{3}}\, {{s}^{4}}}{10871635968}+\frac{{{r}^{2}}\, {{s}^{4}}}{67108864}-\frac{7 {{r}^{10}}\, {{s}^{3}}}{402653184}-\frac{155 {{r}^{9}}\, {{s}^{3}}}{10871635968}-\frac{13241 {{r}^{8}}\, {{s}^{3}}}{97844723712}-\frac{5927 {{r}^{7}}\, {{s}^{3}}}{32614907904}-\frac{5927 {{r}^{6}}\, {{s}^{3}}}{32614907904}-\frac{13241 {{r}^{5}}\, {{s}^{3}}}{97844723712}-\frac{155 {{r}^{4}}\, {{s}^{3}}}{10871635968}-\frac{7 {{r}^{3}}\, {{s}^{3}}}{402653184}+\frac{{{r}^{10}}\, {{s}^{2}}}{67108864}+\frac{19 {{r}^{9}}\, {{s}^{2}}}{1207959552}+\frac{2711 {{r}^{8}}\, {{s}^{2}}}{24461180928}+\frac{1913 {{r}^{7}}\, {{s}^{2}}}{48922361856}+\frac{2711 {{r}^{6}}\, {{s}^{2}}}{24461180928}+\frac{19 {{r}^{5}}\, {{s}^{2}}}{1207959552}+\frac{{{r}^{4}}\, {{s}^{2}}}{67108864}-\frac{{{r}^{10}} s}{134217728}-\frac{163 {{r}^{9}} s}{10871635968}-\frac{539 {{r}^{8}} s}{24461180928}-\frac{539 {{r}^{7}} s}{24461180928}-\frac{163 {{r}^{6}} s}{10871635968}-\frac{{{r}^{5}} s}{134217728}+\frac{{{r}^{10}}}{402653184}+\frac{{{r}^{9}}}{1358954496}+\frac{229 {{r}^{8}}}{48922361856}+\frac{{{r}^{7}}}{1358954496}+\frac{{{r}^{6}}}{402653184}$
\end{flushleft}

The elimination of the ideal generated by $I_9$ and $I_{18}$ with respect to $r$ and $s$ is the 0 ideal, which shows that $I_9$ and $I_{18}$ are algebraically independent.
\end{proposition}

We can compute the invariants of $C_6$ by substitute $s=1-r$ into the invariants of $C_3$.

\begin{proposition} The Dixmier invariants of 
$$C_6(r): y^3=x(x-1)(x-r)(x-1+r)$$
are

$$I_3=I_6=I_{12}=I_{15}=0$$

\begin{flushleft}
$I_9=-\frac{65 {{r}^{8}}-260 {{r}^{7}}+1150 {{r}^{6}}-2540 {{r}^{5}}+3959 {{r}^{4}}-3988 {{r}^{3}}+2326 {{r}^{2}}-712 r+89}{331776}$

\

$I_{18}=\frac{25 {{r}^{16}}}{3057647616}-\frac{25 {{r}^{15}}}{382205952}+\frac{1325 {{r}^{14}}}{3057647616}-\frac{1925 {{r}^{13}}}{1019215872}+\frac{79229 {{r}^{12}}}{12230590464}-\frac{105737 {{r}^{11}}}{6115295232}+\frac{447307 {{r}^{10}}}{12230590464}-\frac{31385 {{r}^{9}}}{509607936}+\frac{998905 {{r}^{8}}}{12230590464}-\frac{57233 {{r}^{7}}}{679477248}+\frac{817465 {{r}^{6}}}{12230590464}-\frac{123275 {{r}^{5}}}{3057647616}+\frac{221939 {{r}^{4}}}{12230590464}-\frac{1337 {{r}^{3}}}{226492416}+\frac{16037 {{r}^{2}}}{12230590464}-\frac{17 r}{95551488}+\frac{17}{1528823808}$
\end{flushleft}

The elimination of the ideal generated by $I_9$ and $I_{18}$ with respect to $r$ is irreducible and generated by

\begin{flushleft}
 $4000000 I_9^{8}-1998092052000 I_9^{7}-676000000 I_9^{6} I_{18}-71509053768117831
      I_9^{6}+224328787434000 I_9^{5} I_{18}+42841500000 I_9^{4}
      I_{18}^{2}-395361312253919627346 I_9^{5}+8460248600243212740 I_9^{4}
      I_{18}-8372335651553250 I_9^{3} I_{18}^{2}-1206702250000 I_9^{2}
      I_{18}^{3}+36392104317997507611465 I_9^{4}+31914880192757153442492 I_9^{3}
      I_{18}-332936970436116610650 I_9^{2} I_{18}^{2}+103850637726127500 I_9
      I_{18}^{3}+12745792515625 I_{18}^{4}-826890695963630262273456
      I_9^{3}-9875439964247275663003440 I_9^{2} I_{18}-644187721569909674246640 I_9
      I_{18}^{2}+4362752394549791982000 I_{18}^{3}-168880832609781468337056
      I_9^{2}+30826420907787244648372032 I_9 I_{18}+474410438868202394564990304
      I_{18}^{2}+2545539129474834804480 I_9+6939213188282316797541120
      I_{18}+960605665900794374400$.
\end{flushleft}

\end{proposition}

For $C_9$, we have

\begin{proposition} The Dixmier invariants of 
$$C_9: y^3=x(x^3-1)$$
are all zero.

\end{proposition}

\section{The Bitangents of $C_3$, $C_6$ and $C_9$}
\subsection{The Bitangents of Plane Quartics}\label{bitan}

The classical theory of plane quartics says that it has 28 bitangents. Recall that a line $L$ is a {\bf bitangent} of a plane curve $C$ if it tangents $C$ at two points $p_1,p_2$ up to multiplicity.

Explicitly, let $f=f(x,y,z)\in K[x,y,z]_4$ be the equation of a plane quartic $C$. 
Let $L: ax+by+cz=0$, $a,b,c\in K$ be a line in ${\mathbb P}^2_{(x,y,z)}$. Thus the point $(a,b,c)\in {\mathbb P}^2_{(a,b,c)}$ determines the line $L$. 
So without lost of generality, we can assume that $c\neq 0$, and say $c=1$. 
This time $L: ax+by+z=0$ gives the condition $z=-ax-by$. 
Substitute this relation into $f(x,y,z)$ we have a quadratic form $f(x,y,-ax-by)\in R[x,y,z]_2$ where $R=K[a,b]$. 
If $L$ is a bitangent for some $a,b\in K$, then there exist $\lambda_0,\lambda_1,\lambda_2\in K$ such that 
\begin{equation}\label{lambda012}f(x,y,-ax-by)=(\lambda_0x^2+\lambda_1xy+\lambda_2y^2)^2.\end{equation}
\begin{definition}\label{IJ}For any quartic $f\in K[x,y,z]_4$, let $I(f)$ be the ideal of $K[a,b,\lambda_0,\lambda_1,\lambda_2]$ generated by comparing the coefficients of both sides of the monomials of  $x,y$ in the expansion of (\ref{lambda012}). 
Let $J(f)$ be elimination ideal of $I$ with respect to $\lambda_0,\lambda_1,\lambda_2$ in $K[a,b]$. \end{definition}  
The ideal $J(f)$ gives the conditions of $L$ being a bitangent of $C$. In general one cannot solve $a,b$ over ${\mathbb Q}$, and even there exists $L$ such that $a,b\in {\mathbb Q}$, the tangency points $p_1,p_2$ are not $\mathbb Q$-rational points of $C$.

There is a description of the relative positions of the bitangents of $C$. Let $L_1,\ldots, L_{28}$ be the bitangents of $C$, be careful that the number 28  counts the overlaps of the bitangents. 
Let $L_i,L_j,L_k$, where $i,j,k=1,\ldots , 28$ are distinct, be a triple of  bitangents. For each $L_\nu$, $\nu=1,\ldots , 28$, let $p_{\nu_1}$, $p_{\nu_2}$ be the two tangency points of $L_\nu$ and $C$. Then $L_i,L_j,L_k$ determine 6 points on $C$. Generically a plane conic is determined by 5 points. 
\begin{definition}If the 6 points $p_{i_1}, p_{i_2}, p_{j_1}, p_{j_2}, p_{k_1}, p_{k_2}$ lie on a plane conic, then we say the triple $L_i, L_j, L_k$ are {\bf sysgetic}, or else we say they are {\bf asysgetic}.\end{definition}

\subsection{The Bitangents of $C_3$, $C_6$ and $C_9$}\label{c369bit}

Before we use the computer to comply the algorithm above, let us observe an obvious bitangent of 
$$C_3(r,s): y^3=x(x-1)(x-r)(x-s).$$

In the algorithm above, we considered the generic case on the affine chart $z\neq 0$. But if we expand $C_3$ and homogenize it with respect to $z$, then we have
\begin{equation}\label{homc3}F_3(r,s): r s x\, {{z}^{3}}-r s\, {{x}^{2}}\, {{z}^{2}}-s\, {{x}^{2}}\, {{z}^{2}}-r\, {{x}^{2}}\, {{z}^{2}}+{{y}^{3}} z+s\, {{x}^{3}} z+r\, {{x}^{3}} z+{{x}^{3}} z-{{x}^{4}}.\end{equation}
Substitute $z=0$ into (\ref{homc3}) we get $x^4=(x^2)^2$, which is a square. Thus $z=0$ is a bitangent of $C_3$. To compute the tangent point, we observe that $x^2=0$ implies that $x=0$. Substitute $x=0,z=0$ into (\ref{homc3}) we get 0. This means that the intersection of $C_3$ and the line $z=0$ is the point $(0,y,0)$, or $(0,1,0)\in {\mathbb P}^2_{(x,y,z)}$. This is the only solution of the point of tangency, so the 2 points of tangency coincide.

Beyond this bitangent, there are another 27 bitangents of $C_3$. Let $J(C_3)$ be the ideal defined as Definition \ref{IJ}. This time the coefficient list becomes $K[r,s]$, but we still can define $J(C_3)$ by the same analogos. We can compute the primary decomposition of $J(C_3)$ using {\tt Macaulay2}. The inputs are as the following.

\begin{verbatim}
R = QQ[r,a,b,k_0,k_1,k_2][x,y,z]
R
f = -r^2*x*z^3+r*x*z^3+r^2*x^2*z^2-r*x^2*z^2-x^2*z^2+y^3*z+2*x^3*z-x^4
g = (k_0*x^2+k_1*x*y+k_2*y^2)^2
h = substitute(f,{z => -a*x-b*y})
H= h-g
Coe = coefficients H
L = flatten entries Coe#1
S = QQ[r,a,b,k_0,k_1,k_2]
I = ideal L
psi=map(S,R)
phi=map(R,S)
J = psi I
E=eliminate(J,{k_0,k_1,k_2})
T = QQ[r,a,b]
xi=map(T,S)
U = xi E
degree U
primaryDecomposition U
\end{verbatim}
The primary decomposition of $J(C_3)$ has two components, one of them is the ideal $<a=0,b=0>$, which gives the bitangent $z=0$. Another component is irreducible in general. Let $J'$ be this component, and let $J'_a$ be the elimination of $J'$ with respect to $b$. Then one can see that $a$ satisfies the degree 9 equation

\begin{flushleft}
$r^{4} s^{4} a^{9}-12 r^{4} s^{3} a^{7}-12 r^{3} s^{4}
      a^{7}-8 r^{4} s^{3} a^{6}-8 r^{3} s^{4} a^{6}-12 r^{3} s^{3} a^{7}-8 r^{4}
      s^{2} a^{6}-120 r^{3} s^{3} a^{6}-8 r^{2} s^{4} a^{6}+30 r^{4} s^{2}
      a^{5}-156 r^{3} s^{3} a^{5}+30 r^{2} s^{4} a^{5}-8 r^{3} s^{2} a^{6}-8
      r^{2} s^{3} a^{6}+48 r^{4} s^{2} a^{4}-96 r^{3} s^{3} a^{4}+48 r^{2} s^{4}
      a^{4}-156 r^{3} s^{2} a^{5}-156 r^{2} s^{3} a^{5}+16 r^{4} s^{2} a^{3}-32
      r^{3} s^{3} a^{3}+16 r^{2} s^{4} a^{3}+48 r^{4} s a^{4}-168 r^{3} s^{2}
      a^{4}-168 r^{2} s^{3} a^{4}+48 r s^{4} a^{4}+30 r^{2} s^{2} a^{5}+68 r^{4}
      s a^{3}-68 r^{3} s^{2} a^{3}-68 r^{2} s^{3} a^{3}+68 r s^{4} a^{3}-96
      r^{3} s a^{4}-168 r^{2} s^{2} a^{4}-96 r s^{3} a^{4}+24 r^{4} s a^{2}-24
      r^{3} s^{2} a^{2}-24 r^{2} s^{3} a^{2}+24 r s^{4} a^{2}+16 r^{4} a^{3}-68
      r^{3} s a^{3}-216 r^{2} s^{2} a^{3}-68 r s^{3} a^{3}+16 s^{4} a^{3}+48
      r^{2} s a^{4}+48 r s^{2} a^{4}+24 r^{4} a^{2}+24 r^{3} s a^{2}-96 r^{2}
      s^{2} a^{2}+24 r s^{3} a^{2}+24 s^{4} a^{2}-32 r^{3} a^{3}-68 r^{2} s
      a^{3}-68 r s^{2} a^{3}-32 s^{3} a^{3}+9 r^{4} a+12 r^{3} s a-42 r^{2}
      s^{2} a+12 r s^{3} a+9 s^{4} a-24 r^{3} a^{2}-96 r^{2} s a^{2}-96 r s^{2}
      a^{2}-24 s^{3} a^{2}+16 r^{2} a^{3}+68 r s a^{3}+16 s^{2} a^{3}+12 r^{3}
      a-12 r^{2} s a-12 r s^{2} a+12 s^{3} a-24 r^{2} a^{2}+24 r s a^{2}-24
      s^{2} a^{2}+8 r^{3}-8 r^{2} s-8 r s^{2}+8 s^{3}-42 r^{2} a-12 r s a-42
      s^{2} a+24 r a^{2}+24 s a^{2}-8 r^{2}+16 r s-8 s^{2}+12 r a+12 s a-8 r-8
      s+9 a+8=0$
\end{flushleft}
which is able to be output by Macaulay2. However, the elimination with respect to $a$ is out of the capability of the processor. This equation is irreducible over ${\mathbb Q}$. In the following cases, we try to find explicit bitangents for special cases of $C_3(r,s)$.

\begin{theorem}
The curve \begin{equation}\label{c9}C_9: y^3=x(x^3-1)\end{equation}
has all 28 explicit symbolic solutions for the bitangents.
\end{theorem}
{\it Proof}\quad Let $J(C_9)$ be the ideal of $K[a,b]$ defined as Definition \ref{IJ}. Let $J'$ be the component of $J(C_9)$ beyond $<a=0,b=0>$. Let $J_a'$ be the elimination of $J'$ with respect to $b$. Then $a$ satisfies the following equation.
\begin{equation}
a^{9}-96 a^{6}+48 a^{3}+64.
\end{equation}
Let $u=a^3$, then $u$ satisfies the cubic equation
\begin{equation}
u^3-96u^2+48u+64.
\end{equation}
This equation is solvable. For example, using Maxima, we have

\begin{align*}
&u_1=-\frac{\left( \sqrt{3} {\rm \, i}+1\right) \, {{\left( 32 \cdot {{3}^{\frac{7}{2}}} {\rm \, i}+31968\right) }^{\frac{2}{3}}}-64 {{\left( 32 \cdot {{3}^{\frac{7}{2}}} {\rm \, i}+31968\right) }^{\frac{1}{3}}}-112\cdot {{3}^{\frac{5}{2}}} {\rm \, i}+1008}{2 {{\left( 32 \cdot {{3}^{\frac{7}{2}}} {\rm \, i}+31968\right) }^{\frac{1}{3}}}},
\\&u_2=\frac{\left( \sqrt{3} {\rm \, i}-1\right) \, {{\left( 32 \cdot {{3}^{\frac{7}{2}}} {\rm \, i}+31968\right) }^{\frac{2}{3}}}+64 {{\left( 32 \cdot {{3}^{\frac{7}{2}}} {\rm \, i}+31968\right) }^{\frac{1}{3}}}-112\cdot {{3}^{\frac{5}{2}}} {\rm \, i}-1008}{2 {{\left( 32 \cdot {{3}^{\frac{7}{2}}} {\rm \, i}+31968\right) }^{\frac{1}{3}}}},
\\&u_3=\frac{{{\left( 32 \cdot {{3}^{\frac{7}{2}}} {\rm \, i}+31968\right) }^{\frac{2}{3}}}+32 \cdot {{\left( 32 \cdot {{3}^{\frac{7}{2}}} {\rm \, i}+31968\right) }^{\frac{1}{3}}}+1008}{{{\left( 32 \cdot {{3}^{\frac{7}{2}}} {\rm \, i}+31968\right) }^{\frac{1}{3}}}}\end{align*}
Taking the cube root of each $u_i$ we can get all 9 solutions of $a$.

Similarly we have an equation
$$b^{27}   - 29496b^{18}   + 401808b^9   - 64=0$$
and let $v=b^9$ we have a cubic equation
$$v^3 - 29496v^2+401808v-64=0.$$
This time one has to take the ninth root of all the three solutions $v_i$'s , $i=1,2,3$ of this equation. At the end, one has to judge which pairs $(a,b)$ among the solutions give a bitangent $ax+by+z=0$ of the original curve. We list the {\tt Macaulay2} input as the following.
\begin{verbatim}
R = QQ[r,s,b,k_0,k_1,k_2][x,y,z]
R
f = r*s*x*z^3-r*s*x^2*z^2-s*x^2*z^2-r*x^2*z^2+y^3*z+s*x^3*z+r*x^3*z+x^3*z-x^4
g = (k_0*x^2+k_1*x*y+k_2*y^2)^2
h = substitute(f,{z => -b*y})
H= h-g
Coe = coefficients H
L = flatten entries Coe#1
S = QQ[r,s,b,k_0,k_1,k_2]
I = ideal L
psi=map(S,R)
phi=map(R,S)
J = psi I
E=eliminate(J,{k_0,k_1,k_2})
T = QQ[r,s,b]
xi=map(T,S)
U = xi E
degree U
primaryDecomposition U
\end{verbatim} \qquad $\blacksquare$

There is no canonical method to find explicit bitangents for special cases. Our observation is that we can try to find $r,s\in K$ such that the bitangent is ``horizontal", that is, for those bitangents such that $a=0$. The equation of the bitangent becomes $bx+z=0$. Repeat the same idea in Section 
\ref{bitan}, we get the following result.

\begin{theorem}\label{c3h}The family $C_3$ has a horizontal bitangent when $r-s=\pm 1$ or $r+s=1$. In each of these cases, the slope $b$ satisfies a cubic equation whose coefficients are polynomials of $s$, thus there are $3$ horizontal bitangents.
\end{theorem}

{\it Proof}\qquad Let $F_3$ be the polynomial defined in (\ref{homc3}). When $a=0$, we have $L:by+z=0$. Then $z=-by$. Using the same idea as in Section \ref{bitan},  
we have the equation
\begin{equation}\label{by}F_3(x,y,-by)=(\lambda_0x^2+\lambda_1xy+\lambda_2y^2)^2.\end{equation}
Let $I(F_3)$ be the ideal of $R[b,\lambda_0,\lambda_1,\lambda_2]$ generated by comparing the coefficients of both sides of the monomials of  $x,y$ in the expansion of (\ref{by}). 
Let $J(F_3)$ be elimination ideal of $I(F_3)$ with respect to $\lambda_0,\lambda_1,\lambda_2$ in $R[b]$. Then the primary decomposition of $J(F_3)$ as an ideal in $K[r,s,b]$
is
\begin{equation}\label{pdc}\begin{split}&\langle b \rangle, \quad  \langle r-s-1,s^{2} b^{3}-4\rangle\\
& \langle r+s-1,s^{4} b^{3}-2 s^{3} b^{3}+s^{2} b^{3}-4\rangle \\
& \langle r-s+1,s^{2} b^{3}-2 s b^{3}+b^{3}-4\rangle\qquad \blacksquare\end{split}\end{equation}
The first ideal of (\ref{pdc}) corresponds to the bitangent $z=0$. The third ideal of (\ref{pdc}) gives $r+s-1=0$, which implies $s=r-1$, this is the family $C_6$. Furthermore, we have a result on the positions of the horizontal bitangents of $C_6$.

\begin{theorem}\label{c6bitan}
The three horizontal bitangents of $C_6$ form an asyzygetic triple.
\end{theorem}

{\it Proof}\qquad Let
$$F_6: -{{r}^{2}} x\, {{z}^{3}}+r x\, {{z}^{3}}+{{r}^{2}}\, {{x}^{2}}\, {{z}^{2}}-r\, {{x}^{2}}\, {{z}^{2}}-{{x}^{2}}\, {{z}^{2}}+{{y}^{3}} z+2 {{x}^{3}} z-{{x}^{4}}\in R[x,y,z]_4$$
be the homogenization of $C_6$ with respect to $z$ where $R=K[r]$. As before, we have the equation
\begin{equation}\label{f6}F_6(x,y,-by)=(\lambda_0x^2+\lambda_1xy+\lambda_2y^2)^2\end{equation}
Let $I(F_6)$ be the ideal of $R[b,\lambda_0,\lambda_1,\lambda_2]$ generated by comparing the coefficients of both sides of the monomials of  $x,y$ in the expansion of (\ref{f6}).
In Theorem \ref{c3h} we have proved that for $C_6$ the condition of being a horizontal bitangent for the line $bx+z=0$ is given by the ideal
$$\langle\, r+s-1,\, s^{4} b^{3}-2 s^{3} b^{3}+s^{2} b^{3}-4\, \rangle.$$
Substitute $s=1-r$ into the second generator of this ideal, we have a relation
$$p(r,b)= {{b}^{3}}\, {{r}^{4}}-2 {{b}^{3}}\, {{r}^{3}}+{{b}^{3}}\, {{r}^{2}}-4 $$

This time, let ${\mathscr J}(F_6)$ be intersection of the elimination ideal of $I(F_6)$ with respect to $r,b$ in $K[\lambda_0,\lambda_1,\lambda_2]$ and the ideal $\langle p(r,b)\rangle$. 
Macaulay2 outputs $${\mathscr J}(F_6)=\langle\  \rangle,$$
which means that generically there is no conic $\lambda_0x^2+\lambda_1xy+\lambda_2y^2$ satisfies the conditions of passing through the 6 tangent points at the same time.\qquad $\blacksquare$

\begin{remark}In general, there is another way to check whether 6 points lie on a common conic in ${\mathbb P}^2$. Let $p_i=(x_i,y_i,z_i)\in {\mathbb P}^2_{(x,y,z)}, i= 1,\ldots ,6$ be 6 points in the projective plane. Let ${\bf V}$ be the Veronese map
$$\begin{array}{cccc}
{\bf V}: & {\mathbb P}^2_{(x,y,z)} & \longrightarrow & {\mathbb P}^5 \\
& (x,y,z) & \longmapsto & (x^2,y^2,z^2,xy,yz,zx)
\end{array}.$$ If we regard ${\bf V}(p)$ as a row matrix for any $p=(x,y,z)\in {\mathbb P}^2$, then for the given 6 poins $p_1,...,p_6$, we have a $6\times 6$ matrix $$V:=\begin{pmatrix}
{\bf V}(p_1) \\
{\bf V}(p_2) \\
{\bf V}(p_3) \\
{\bf V}(p_4) \\
{\bf V}(p_5) \\
{\bf V}(p_6) 
\end{pmatrix} = \begin{pmatrix}
x_1^2 & y_1^2 & z_1^2 & x_1y_1 & y_1z_1 & z_1x_1 \\
x_2^2 & y_2^2 & z_2^2 & x_2y_2 & y_2z_2 & z_2x_2 \\
x_3^2 & y_3^2 & z_3^2 & x_3y_3 & y_3z_3 & z_3x_3 \\
x_4^2 & y_4^2 & z_4^2 & x_4y_4 & y_4z_4 & z_4x_4 \\
x_5^2 & y_5^2 & z_5^2 & x_5y_5 & y_5z_5 & z_5x_5 \\
x_6^2 & y_6^2 & z_6^2 & x_6y_6 & y_6z_6 & z_6x_6  
\end{pmatrix}.$$

For our problem, let $p_1,...,p_6$ be the 6 points of tangency of the three horizontal bitangents in Theorem \ref{c6bitan}. From the proof of Theorem \ref{c6bitan} we see that there is a symbolic solution of these three bitangents, and since the algorithm of finding the points of tangency is essentially solving a quadratic equation, we can find the symbolic solutions of the points of tangency. But this algorithm costs too much for a popular processor. We can compute it in special values. For example, let $r=\frac{1}{8}$, we can compute the  determinant using Maxima, the result is\footnote{This value could be simplified, we put the original result from Maxima.}
$$V=-\frac{\sqrt{-25 \sqrt{3} {\rm \, i}-25}\, \sqrt{25 \sqrt{3} {\rm \, i}-25}\, \left( 12005 {{2}^{\frac{10}{3}}}\, {{3}^{\frac{7}{2}}}\, {{4}^{\frac{2}{3}}} {\rm \, i}+324135 {{2}^{\frac{10}{3}}}\, {{4}^{\frac{2}{3}}}\right) }{{{2}^{\frac{247}{6}}}\, \sqrt{3} {\rm \, i}+{{2}^{\frac{247}{6}}}}$$
which is not zero.

\end{remark}

\section{Discussion on the Matrix Representation Problem}

We discuss the matrix representation problem of the curves $C_3$ and $C_6$ using the idea in \cite{Sturmfels2}. In order to coincide the notations with respect to \cite{Sturmfels2}, we exchange $y$ and $z$, and write $C_3$ as
$$C_3:\quad z^3=x(x-1)(x-r)(x-s).$$
Homogenize $C_3$ with respect to $y$ we have
\begin{equation}\label{C3new}C_3:\quad f(x,y,z):=F_3(r,s)=x(x-y)(x-ry)(x-sy)-yz^3=0.\end{equation}
This time we have
\begin{equation}\label{condition}f(x,0,0)=x^4\quad \mbox{and}\quad f(x,y,0)=\prod_{i=1}^4 (x+\beta_iy)\end{equation}
where $\beta_1=0,\beta_2=-1,\beta_3=-r,\beta_4=-s$.
The matrix representation problem for $C_3$ asks whether the polynomial $f(x,y,z)$ in ($\ref{C3new}$) could be written of the form
$$f(x,y,z)={\rm det}(xA+yB+zC)$$
where $A,B,C$ are symmetric matrices.
 Here the entries of the matrices $A,B$ and $C$ belong to the algebraic closure of the rational function field $K(r,s)$. According to Section 2 in  \cite{Sturmfels2}, if (\ref{condition}) holds, then one can assume that
$$
 A =\begin{pmatrix}1 &&& \\
 &1 && \\ && 1 & \\ &&& 1 \end{pmatrix}, \quad B = \begin{pmatrix}0 &&& \\
 & -1 && \\ && -r & \\ &&& -s \end{pmatrix} 
,\quad C=\begin{pmatrix}{c_{11}} & {c_{12}} & {c_{13}} & {c_{14}}\\
{c_{12}} & {c_{22}} & {c_{23}} & {c_{24}}\\
{c_{13}} & {c_{23}} & {c_{33}} & {c_{34}}\\
{c_{14}} & {c_{24}} & {c_{34}} & {c_{44}}\end{pmatrix}.$$
and we also have that
\begin{equation}\label{cii}c_{ii}=\beta_i\cdot \frac{\frac{\partial f}{\partial z}(-\beta_i,1,0)}{\frac{\partial f}{\partial y}(-\beta_i,1,0)},\quad i=1,2,3,4.\end{equation}
But for (\ref{C3new}) we have $\frac{\partial f}{\partial z}=-3yz^2,$
which implies if $z=0$, then $c_{ii}=0$ for $i=1,2,3,4$ by (\ref{cii}).

For convinience we denote
$$D=\begin{pmatrix} & {c_{12}} & {c_{13}} & {c_{14}}\\
& & {c_{23}} & {c_{24}}\\
&& & {c_{34}}\\
&&& \end{pmatrix} = \begin{pmatrix} & a & b & d\\
& & c & e\\
&& & f\\
&&& \end{pmatrix},$$
then $C=D+\,^{\sf t}D$ where $\,^{\sf t}D$ is the matrix transpose of $D$ since $c_{ii}=0$ for $i=1,2,3,4$.

Using Maxima, we directly compute the coefficients of $$\det (xA+yB+zC)=\det \begin{pmatrix} x & a z & b z & d z\\
a z & x-y & c z & e z\\
b z & c z & x-r y & f z\\
d z & e z & f z & x-s y\end{pmatrix}$$ and compare the coefficients with $f(x,y,z)$ in (\ref{C3new}), the output is a system of equations
\begin{align}
& -c^2s-b^2s-a^2s-e^2r-d^2r-a^2r-f^2-d^2-b^2=0, \label{oeq1} \\
& {{a}^{2}} r s+{{b}^{2}} s+{{d}^{2}} r =0,\label{oeq2}\\ 
&2 a b c s+2 a d e r+2 b d f-1=0,\label{oeq3}\\
& {{f}^{2}}+{{e}^{2}}+{{d}^{2}}+{{c}^{2}}+{{b}^{2}}+{{a}^{2}}=0, \label{oeq4}\\
& -2 c e f-2 b d f-2 a d e-2 a b c=0,\label{oeq5}\\
& -{{a}^{2}}\, {{f}^{2}}+2 a b e f+2 a c d f-{{b}^{2}}\, {{e}^{2}}+2 b c d e-{{c}^{2}}\, {{d}^{2}}=0\label{oeq6}.
\end{align} 
We add the first equation with the fourth one, and rewrite the system of as 6 equations
\begin{align}
& {{a}^{2}} r s+{{b}^{2}} s+{{d}^{2}} r =0, \label{e1}\\ 
& a^2(1-r)(s-1-s)+c^2(1-s)+e^2(1-r)=0 \label{e2}\\
&2 a b c s+2 a d e r+2 b d f-1=0,\label{e3}\\
& {{f}^{2}}+{{e}^{2}}+{{d}^{2}}+{{c}^{2}}+{{b}^{2}}+{{a}^{2}}=0,\label{e4}\\
&  c e f+ b d f+ a d e+ a b c=0,\label{e5}\\
& a^2f^2-2af(be+cd)+(be-cd)^2=0\label{e6}.
\end{align}
of the 6 variables $a,b,c,d,e,f$.

It is too complicated to solve this entire system. Our computation are proceeded under the following principle:

\begin{itemize}
\item We only seek for one solution to the equation system (\ref{e1})-(\ref{e6}), thus if there is an "either-or" argument in any step, we can choose one of them as our solution.
\end{itemize}
We eliminate $a,f$, and get a system of 4 equations with respect to the 4 variables $b,c,d,e$.

\begin{proposition}The equation system 
\begin{align}
&\frac{b^2}{r}+\frac{c^2}{r-1}+\frac{d^2}{s}+\frac{e^2}{s-1}=0, \label{eb1} \\
&\frac{(b-e)^4}{be}=\frac{(c+d)^4}{cd}, \label{eb2}\\
&(bc+de)(bd+ce)=\left(2(\sqrt{be}+\sqrt{cd})\cdot\begin{vmatrix}
bcs+der & bd \\ bc(1-s)+de(1-r) & ce
\end{vmatrix}\right)^2,\label{eb3}\\
&\frac{(bd+ce)^2+(bc+de)^2}{\ \begin{vmatrix}
bcs+der & bd \\ bc(1-s)+de(1-r) & ce
\end{vmatrix}^2\ }+(b^2+c^2+d^2+e^2)=0\label{eb4}
\end{align}
with respect to the variables $b,c,d,e$ give answers to the equation system (\ref{e1})-(\ref{e6}) where 
\begin{equation}\label{eb5}a= \frac{bd+ce}{\ \begin{vmatrix}
bcs+der & bd \\ bc(1-s)+de(1-r) & ce
\end{vmatrix}\ },\quad f=-\frac{bc+de}{\ \begin{vmatrix}
bcs+der & bd \\ bc(1-s)+de(1-r) & ce
\end{vmatrix}\ }.\end{equation}
\end{proposition}

{\it Proof}\qquad First, the equation (\ref{eb1}) is simply from $\frac{1}{rs}$(\ref{e1})$-\frac{1}{(1-r)(1-s)}$(\ref{e2}). 

Next we regard $a,f$ as unknowns and $b,c,d,e,r,s$ as constants. The solution (\ref{eb5}) is the answer to the linear system (\ref{e3}) and (\ref{e5}). From (\ref{e5}) we also have
\begin{equation}\label{es1}\frac{a}{f}=-\frac{bd+ce}{bc+de},\quad \frac{f}{a}=-\frac{bc+de}{bd+ce1}.
\end{equation}
Substitude (\ref{es1}) into (\ref{e4}) we have ({\ref{eb4}}).

Let $g=af$, then (\ref{e6}) becomes a quadratic equation
$$g^2-2(be+cd)g+(be-cd)^2=0$$
of $g$ whose solution is
$$
af=(\sqrt{be}\pm \sqrt{cd})^2 $$
As before, for ``$\pm$" we choose $+$, which is
\begin{equation}\label{eaf}af=(\sqrt{be}\pm \sqrt{cd})^2
\end{equation}
Substitude (\ref{eb5}) into (\ref{eaf}) we get (\ref{eb3}).

Last, let us prove (\ref{eb2}). The quadratic equation (\ref{e4}) and the linear equation (\ref{e5}) have an answer
\begin{equation}\label{eaf2}\begin{split}
&a=-\frac{\sqrt{-1} (c e+bd)\, \sqrt{{{e}^{2}}+{{d}^{2}}+{{c}^{2}}+{{b}^{2}}}\, \sqrt{\left( {{d}^{2}}+{{c}^{2}}\right) \, {{e}^{2}}+4 b c d e+{{b}^{2}}\, {{d}^{2}}+{{b}^{2}}\, {{c}^{2}}}}{\left( {{d}^{2}}+{{c}^{2}}\right) \, {{e}^{2}}+4 b c d e+{{b}^{2}}\, {{d}^{2}}+{{b}^{2}}\, {{c}^{2}}}, \\ &f=\frac{\sqrt{-1} (bc+d e)\, \sqrt{{{e}^{2}}+{{d}^{2}}+{{c}^{2}}+{{b}^{2}}}\, \sqrt{\left( {{d}^{2}}+{{c}^{2}}\right) \, {{e}^{2}}+4 b c d e+{{b}^{2}}\, {{d}^{2}}+{{b}^{2}}\, {{c}^{2}}}}{\left( {{d}^{2}}+{{c}^{2}}\right) \, {{e}^{2}}+4 b c d e+{{b}^{2}}\, {{d}^{2}}+{{b}^{2}}\, {{c}^{2}}}.\end{split}
\end{equation}
From (\ref{eaf2}) we have
$$
af=\frac{P}{\left( {{d}^{2}}+{{c}^{2}}\right) \, {{e}^{2}}+4 b c d e+{{b}^{2}}\, {{d}^{2}}+{{b}^{2}}\, {{c}^{2}}}
$$
where the numerator $P$ equals to
\begin{align*}&-\left( c d\, {{e}^{4}}-4 b c d\, {{e}^{3}}+6 {{b}^{2}} c d\, {{e}^{2}}-b\, {{d}^{4}} e-4 b c\, {{d}^{3}} e-6 b\, {{c}^{2}}\, {{d}^{2}} e-4 b\, {{c}^{3}} d e-4 {{b}^{3}} c d e-b\, {{c}^{4}} e+{{b}^{4}} c d\right) \\&\, \left( c d\, {{e}^{4}}+4 b c d\, {{e}^{3}}+6 {{b}^{2}} c d\, {{e}^{2}}-b\, {{d}^{4}} e+4 b c\, {{d}^{3}} e-6 b\, {{c}^{2}}\, {{d}^{2}} e+4 b\, {{c}^{3}} d e+4 {{b}^{3}} c d e-b\, {{c}^{4}} e+{{b}^{4}} c d\right), \end{align*}
we take the first factor
\begin{equation}\label{bcd}
 c d\, {{e}^{4}}-4 b c d\, {{e}^{3}}+6 {{b}^{2}} c d\, {{e}^{2}}-b\, {{d}^{4}} e-4 b c\, {{d}^{3}} e-6 b\, {{c}^{2}}\, {{d}^{2}} e-4 b\, {{c}^{3}} d e-4 {{b}^{3}} c d e-b\, {{c}^{4}} e+{{b}^{4}} c d
\end{equation}
Divide ({\ref{bcd}) by $bcde$ and regroup the terms, we prove (\ref{eb2}). \qquad $\blacksquare$

As we reminded, it is hard to continue solving this equation system. Our observation is that for (\ref{eb2}), we have an obvious solution 
\begin{equation}\label{bcde}
e=b,\quad  \mbox{and}\quad d=-c.
\end{equation}
From (\ref{oeq4}) we have 
$b^2+d^2+f^2=-a^2-c^2-e^2$, thus we can rewrite (\ref{oeq1}) as
$$(a^2+b^2+d^2)r+(a^2+b^2+c^2)s=a^2+c^2+e^2.$$
Substitude (\ref{bcde}) into this equation we have
$$r+s=1$$
which means the curve $C_3$ becomes $C_6$ in this situation.

Next, we subsitude (\ref{bcde}) into the equation system (\ref{e1})-(\ref{e6}), then (\ref{e5}) is trivial, and (\ref{e1}) is the same as (\ref{e2}). We have a system of 4 equations
\begin{align}
&a^2rs+b^2s+c^2r=0 \label{ebe1}\\
&2abc(s-r)-2bcf-1=0 \label{ebe2} \\
&a^2+f^2+2(b^2+c^2)=0 \label{ebe3} \\
&a^2f^2-2af(b^2-c^2)+(b^2+c^2)^2=0 \label{ebe4}
\end{align}
of the 4 variables $a,b,c,f$.

\begin{theorem}\label{5.1}
The matrix representation problem for $C_6$ has a symbolic solution as a solution of a degree 6 polynomial of 1 variable whose coefficients are defined over $\overline{K(r,s)}$.
\end{theorem}

{\it Proof}\qquad From (\ref{ebe2}) we have
$$(a(s-r)-f)=\frac{1}{2bc},$$
thus we have
\begin{equation}\label{abcf6}a^2(s-r)^2-2af(s-r)+f^2=\frac{1}{4b^2c^2}.\end{equation}
From (\ref{ebe3}) we have $f^2=-2(b^2+c^2)-a^2$, substitude it into (\ref{abcf6}) we have
\begin{equation}\label{abcf62}a^2[(s-r)^2-1]-2af(s-r)-2(b^2+c^2)=\frac{1}{4b^2c^2}.\end{equation}
From (\ref{ebe1}) we have \begin{equation}\label{a2}a^2=-\frac{b^2}{r}-\frac{c^2}{s}\end{equation}
and from (\ref{ebe4}) we have \begin{equation}\label{af}af=(b+\sqrt{-1}c)^2\end{equation}
 if we take one of the solutions of the quadratic equation with respect to $af$. Substitude them into (\ref{abcf62}), we have
\begin{equation}
4(b^2s+c^2r)-2(s-r)(b+\sqrt{-1}c)^2-2(b^2+c^2)=\frac{1}{4b^2c^2}.
\end{equation}
This is a degree 6 equation with respect to $b$ and $c$. Thus, if we know $q=b/c$, then the theorem is proved.
From (\ref{ebe3}) and (\ref{a2}) we can solve
\begin{equation}\label{f2}f^2=-2(b^2+c^2)+\frac{b^2}{r}+\frac{c^2}{s}.
\end{equation}
The trivial equation
$$(af)^2=a^2\cdot f^2$$
implies that (\ref{af})$^2=$(\ref{a2})$\cdot$(\ref{f2}), which is
\begin{equation}
(b+\sqrt{-1}c)^4=\left(-\frac{b^2}{r}-\frac{c^2}{s}\right)\cdot \left(-2(b^2+c^2)+\frac{b^2}{r}+\frac{c^2}{s}\right)
\end{equation}
This equation is homogeneous of degree 4 with respect to $b$ and $c$, thus if we set $q=b/c$, it will become a degree 4 equation of $q$, which is solvable. \qquad $\blacksquare$

\begin{flushright}

{\sc School of Mathematics

Sun Yat-Sen University

Guangzhou China

510275
}

{\tt liangdun@mail.sysu.edu.cn}
\end{flushright}

\end{document}